\documentclass[titlepage,11pt]{article}
\usepackage{latexsym}
\usepackage{amsfonts}
\usepackage{amsmath}
\newcommand\blackslug{\hbox{\hskip 1pt \vrule width 4pt height 8pt depth 1.5pt
        \hskip 1pt}}
\newcommand\bbox{\hfill \quad \blackslug \bigbreak}

\def\LL{,\ldots,}
\def\cupcup{\cup\cdots\cup}

\title{On a problem of El-Zahar and Erd\H{o}s}
\author{
Tung Nguyen\thanks{Supported by AFOSR grant
A9550-19-1-0187.}\\
Princeton University, Princeton, NJ 08544, USA
\\
\\
Alex Scott\thanks{Research supported by EPSRC grant EP/X013642/1.}\\
Mathematical Institute, University of Oxford, Oxford OX2 6GG, UK
\\
\\
Paul Seymour\thanks{Supported by AFOSR grant
FA9550-22-1-0234, and NSF grant DMS-2154169.}\\
Princeton University, Princeton, NJ 08544, USA
}

\date{February 6, 2023; revised \today}

\newtheorem{thm}{}[section]

\newcommand{\Proof}{\noindent{\bf Proof.}\ \ }

\begin{document}
\maketitle
\begin{abstract}
Two subgraphs $A,B$ of a graph $G$ are {\em anticomplete} if they are vertex-disjoint and there are no edges joining them.
Is it true that if $G$ is a graph with bounded clique number, and sufficiently large chromatic number, then it has two anticomplete 
subgraphs, 
both with large chromatic number? This is a question raised by El-Zahar and  Erd\H{o}s in 1986, and 
remains open. If so, then at least there should be two anticomplete subgraphs both with large minimum degree, and that is one
of our results.  

We prove two variants of this. First, a strengthening: we can ask for one of the two subgraphs to have large 
chromatic number: that is, 
for all $t, c\ge 1$
there exists $d\ge 1$ such that if $G$ has chromatic number at least $d$, and does not contain the complete graph $K_t$ as a subgraph,
then there are
anticomplete subgraphs $A,B$, where $A$ has minimum degree at least $c$ and $B$ has chromatic number at least $c$.

Second, we look at what happens if we replace the hypothesis that $G$ has sufficiently large chromatic number with the hypothesis that $G$ 
has sufficently large minimum degree. This, together with excluding $K_t$, is {\em not} enough to guarantee 
two anticomplete subgraphs both with large minimum degree; but it works if instead of excluding $K_t$ we exclude the complete 
bipartite graph $K_{t,t}$. More exactly:
for all $t, c\ge 1$
there exists $d\ge 1$ such that if $G$ has minimum degree at least $d$, and does not contain the complete bipartite graph 
$K_{t,t}$ as a subgraph, then there are
two anticomplete subgraphs both with minimum degree at least $c$.

\end{abstract}

\section{Introduction}

We begin with some notation. If $G$ is a graph and $A\subseteq V(G)$, $G[A]$ denotes the subgraph induced on $A$. The chromatic number
of $G$ is denoted by $\chi(G)$, the size of its largest clique is denotes by $\omega(G)$, and 
and if $A\subseteq V(G)$, we sometimes write $\chi(A)$ for $\chi(G[A])$.  
If $A,B$ are subsets of $V(G)$, they are {\em anticomplete} if $A\cap B=\emptyset$ and there
are no edges of $G$ between $A$ and $B$. 

There is a well-known problem of El-Zahar and  Erd\H{o}s~\cite{elzahar, erdos}:
\begin{thm}\label{erdosconj}
{\bf Problem:} Is the following true? For all integers $t, c\ge 1$,
there exists $d\ge 1$, such that if $\chi(G)\ge d$ and $\omega(G)< t$,
then there are anticomplete subsets $A,B\subseteq V(G)$ with $\chi(A), \chi(B)\ge c$.
\end{thm}
This remains open. El-Zahar and  Erd\H{o}s proved that under the same hypotheses,  there are anticomplete subsets 
$A,B\subseteq V(G)$ with $\chi(A)\ge 3$ and $\chi(B)\ge c$, but there has been little further progress. (See~\cite{milner}
for results on an analogous question with infinite graphs and infinite chromatic number.)
We remark that if we omit the hypothesis about $\omega(G)$, the result is no longer true, and a large complete graph is a 
counterexample.

Minimal graphs with large 
chromatic number have large minimum degree, and so if \ref{erdosconj} is true, under the same hypotheses there should at least be
anticomplete subsets $A,B\subseteq V(G)$ such that $G[A]$, $G[B]$ have minimum degree at least $c$. This is true, and can be strengthened:
we can require that one of $G[A],G[B]$ has chromatic number at least $c$. We will prove:
\begin{thm}\label{chimainthm}
For all integers $t, c\ge 1$,
there exists $d\ge 1$, such that if $\chi(G)\ge  d$ and $\omega(G)< t$,
then there are anticomplete subsets $A,B\subseteq V(G)$ where $G[A]$ has minimum degree at least $c$ and $\chi(B)\ge c$.
\end{thm}

What if we relax the hypothesis that $\chi(G)$ is large, and just assume that $G$ has large minimum degree? With $\omega(G)$
bounded, can we still
necessarily find anticomplete subsets $A,B\subseteq V(G)$ such that $G[A]$, $G[B]$ have minimum degree at least $c$? No: a large complete
bipartite graph $G$ is a counterexample. For this question, it becomes natural to bound $\tau(G)$
rather than $\omega(G)$, where $\tau(G)$ is the largest integer $t$ such that $G$ contains $K_{t,t}$ as a subgraph.
We will prove:
\begin{thm}\label{mindegmainthm}
For all integers $t, c\ge 1$,
there exists $d\ge 1$, such that if $G$ has minimum degree at least $d$ and $\tau(G)< t$,
then there are anticomplete subsets $A,B\subseteq V(G)$ where $G[A],G[B]$ both have minimum degree at least $c$.
\end{thm}
Finally, we will examine a possible extension of \ref{erdosconj} to tournaments.

\section{Some lemmas}

We denote the number of vertices of a graph $G$ by $|G|$; and
let us say the {\em denseness} of a non-null graph $G$ is $|E(G)|/|G|$. (In some papers this is called ``density'', but density is also frequently used
to mean something else, so we prefer a different word.) The denseness of the null graph is zero. Also, we define the 
{\em minimum degree} of the null graph to be zero.

The next result is well-known and standard.
\begin{thm}\label{mintoav}
Let $d> 0$. Every graph of minimum degree at least $d$ has denseness at least $d/2$; and
every graph of denseness at least $d$ has a subgraph with minimum degree at least $d$.
\end{thm}
\Proof
The first statement is trivial. For the second, let $G$ be a graph with denseness at least $d$, and choose $G$ minimal 
with these properties. Thus $|E(G)|\ge d|G|$. 
If some vertex $v\in V(G)$ has degree at most $d$, then $|G|\ge 2$ (since $G$ has denseness 
at least $d$ and $d>0$), so the graph $G'$ obtained by deleting $v$ is non-null and satisfies
$$|E(G')|\ge |E(G)|-d\ge d|G|-d=d|G'|,$$
contrary to the minimality of $G$. This proves \ref{mintoav}.~\bbox

In view of \ref{mintoav}, we can replace the conditions about minimum degree
in \ref{chimainthm} and \ref{mindegmainthm} with conditions about denseness, and this is a little more convenient.

If $p\ge 1$ is an integer, let us say a {\em $p$-rock} of a graph $G$ is a set $A\subseteq V(G)$ such that
\begin{itemize}
\item $A\ne \emptyset$ and $|E(G[A])|\ge p|A|$
%p\ge c
\item subject to the above, $|A|$ is minimum; and
\item subject to the two conditions above, $|E(G[A])|$ is maximum.
\end{itemize}

We will need:
\begin{thm}\label{critical}
Let $p\ge 1$ be an integer, let $G$ be a graph, and let $A$ be a $p$-rock of $G$.
Then every vertex $v\in V(G)$ not in $A$ has at most $2p+1$ neighbours in $A$.
\end{thm}
\Proof Since $A\ne \emptyset$ and $|E(G[A])|\ge p|A|$, it follows that $G[A]$ has a vertex of degree at least $2p$, and so 
$|A|\ge 2p+1\ge 2$.
Let $v$ have $t$ neighbours in $A$. For $u\in A$, $(A\setminus \{u\})\cup \{v\}$ has the same cardinality as $A$, so it cannot
induce a subgraph with more edges than $G[A]$. Consequently, every vertex of $A$ has at least $t-1$ neighbours in $A$, and indeed, 
at least $t$
such neighbours, unless it is adjacent to $v$. So 
$$|E(G[A])|\ge \left(t(t-1)+ (a-t)t\right)/2= (a-1)t/2,$$
where $a:=|A|$. 

Let $u\in A$;
then from the minimality of $|A|$, $G[A\setminus \{u\}]$ has at most $p(a-1)-1$ edges.
Summing over all $u\in A$, we deduce that
$$(a-2)|E(G[A])| \le  pa(a-1)-a<(pa-1)(a-1).$$
Substituting, it follows that
$$(a-2)(a-1)t/2 <  (pa-1)(a-1),$$
and so $(a-2)t<2(pa-1)$. Since $a\ge 2p+1$ (because $G[A]$ has denseness at least $p$), it follows that
$(a-2)(2p+2)\ge 2(pa-1)$, and so $t<2p+2$. Hence $t\le 2p+1$. This proves \ref{critical}.~\bbox

We remark that the bound of \ref{critical} is tight, because for instance $A$ might be a clique with $2p+1$ vertices.
Our third lemma is rather obvious, but we will use it twice, so we might as well state it explicitly:
\begin{thm}\label{getmatching}
Let $H$ be a graph and $q\ge 1$ an integer. Then there is a partition of $E(H)$ into sets $M_0\LL M_n$ for some $n\ge 0$,
such that 
\begin{itemize}
\item there is a subset $X\subseteq V(H)$ with $|X|\le 2q-2$
such that every edge in $M_0$ is incident with a vertex in $X$; and
\item $M_1\LL M_n$ are all matchings, each with cardinality $q$.
\end{itemize}
\end{thm}
\Proof We use induction on $|E(H)|$. Suppose first that $H$ has no matching with cardinality $q$. Let $M$ a maximal matching of $H$;
then every edge of $H$ has an end in $X$, where $X$ is the set of vertices incident with an edge of $M$, from the maximality of $M$.
Since $|M|\le q-1$ and hence $|X|\le 2q-2$, we may set $M_0=E(H)$ and $n=0$, and the theorem holds. So we may assume that
$H$ has a matching $M$ of cardinality $q$; but then the result follows from the inductive hypothesis applied to the graph
obtained from $H$ by deleting the edges of $M$. This proves \ref{getmatching}.~\bbox

Fourth, we need:
\begin{thm}\label{probability}
Let $H$ be a graph, and let $Z\subseteq V(H)$ such that each vertex of $H$ belongs to $Z$ independently with probability $1/2$.
\begin{itemize}
\item If $M$ is a matching of $H$, the probability that at most $|M|/8$ edges in $M$ have both ends in $Z$
is at most $e^{-|M|/32}$.
\item Let $d\ge 0$, and for each $v\in V(H)$, let $0\le d_v\le d$, and let $m:=\sum_{v\in V(H)}d_v$; then the probability that 
$\sum_{v\in Z}d_v\le m/4$
is at most $e^{-m/(8d)}$.
\end{itemize}
\end{thm}
\Proof The first statement is immediate from Hoeffding's inequality, since each edge of $M$ has both ends in $Z$ independently 
with probability $1/4$. For the second statement, since $D:=\sum_{v\in Z}d_v$ is a sum of independent bounded random variables, and the expected
value of $D$ is $m/2$, we can apply
Hoeffding's inequality, and deduce that the probability that $|D|\le m/4$ is at most 
$$\exp\left(\frac{-m^2}{8\sum_{v\in V(H)} d_v^2}\right).$$
But  $\sum_{v\in V(H)} d_v^2$ is at most $md$, since $\sum_{v\in V(H)} d_v=m$ and each $d_v\le d$; so 
the probability that $|D|\le m/4$ is at most
$e^{-m/(8d)}$. This proves the second statement, and so proves \ref{probability}.~\bbox

\section{The main proofs}

First we prove \ref{chimainthm}, which we restate (in terms of denseness rather than minimum degree, which by \ref{mintoav} 
is equivalent):

\begin{thm}\label{chimainthm2}
For all integers $t, c\ge 1$,
there exists $d\ge 1$, such that if $G$ is a graph with $\chi(G)\ge  d$ and $\omega(G)< t$,
then there are anticomplete subsets $A,B\subseteq V(G)$ where $G[A]$ has denseness at least $c$ and $\chi(B)\ge c$.
\end{thm}
%p integer
%p=32c$
%r=p

%(%e^{-q/32}2^{2r}\le 1/2)
%q>32 log 2(2r+1)
%$2qe^{-q/8}2^{2r}<1$

%(%3|E(G[A])|\ge 16xq)
%(%3 p|A|/2\ge 16xq)
%(3p|A|/2 \ge 32q^2)

%$d>64q^2d'/(3p)%

%$d>p+1+2qd'+2^{24}c 
%$d\ge 64q^2d'/(3p)+c% defn of d

\Proof
We proceed by induction on $t$. If $t\le 2$ we may take $d=2$, because $\chi(G)\le 1$ for every graph $G$ with $\omega(G)\le 1$.
Thus we may assume that $t\ge 3$, and the result holds for $t-1$. Choose $d'\ge 1$ such that
if for every graph $G$, if $\chi(G)\ge d'$ and $\omega(G)< t-1$,
then there are anticomplete subsets $A,B\subseteq V(G)$ where $G[A]$ has denseness at least $c$ and $\chi(B)\ge c$.

Let $p=32c$; choose $q$ such that 
$e^{-q/32} < 2^{-4p-3}$,
and choose $d$ such that 
$$d>\max\left(2p+1+2qd'+2^{2p}c, 8q^2d'/p+c\right).$$
We will show that $d$ satisfies the theorem.

Let $G$ be a graph with $\omega(G)\le t$, such that there do not exist anticomplete subsets $A,B\subseteq V(G)$ 
where $G[A]$ has denseness at least $c$ and $\chi(B)\ge c$. We will prove that $\chi(G)<d$. 
From the inductive hypothesis, it follows that for every vertex
$v$, its set of neighbours $N$ satisfies $\chi(N)\le d'$.
We may assume that $G$ has a non-null subgraph with minimum degree at least $d-1$, because otherwise $\chi(G)<d$ as required.
Since $p\le (d-1)/2$, 
there is a $p$-rock $A$ of $G$.
Let $F:=E(G[A])$.
By \ref{getmatching},
$F$ may be partitioned into $M_0,M_1\LL M_n$ for some $n\ge 0$, such that
\begin{itemize}
\item there exists $X\subseteq A$ with $|A|\le 2q-2$ such that every edge in $M_0$ has an end in $X$; and
\item $M_1\LL M_n$ are all matchings of cardinality $q$.
\end{itemize}
We may assume that:
\\
\\
\noindent(1) {\em $|A|\ge 8q^2/p$.}
\\
\\
Suppose not. Then the set of vertices of $G$ with a neighbour in $A$ (this set includes $A$, from the minimality of $A$) 
has chromatic number at most $d'|A|\le  8d'q^2/p$;
and the set with no neighbour in $A$ (and that therefore do not belong to $A$) has chromatic number less than $c$, since it is 
anticomplete to $A$ and $p\ge c$. Thus $\chi(G)< 8d'q^2/p+c\le d$ as required. This proves (1).
%64d'q^2/(3p)+c\le d$

\bigskip

Let $\mathcal{I}$ be the set of all subsets of $\{1\LL 4p+2\}$ with cardinality $2p+1$.
\\
\\
(2) {\em There is a partition of $A\setminus X$ into $4p+2$ subsets $A_1\LL A_{4p+2}$, such that for each $I\in \mathcal{I}$,
at least
$|F|/32\ge p|A|/32$ edges have both ends in $X\cup \bigcup_{i\in I}A_i$.}
\\
\\
For each $v\in A\setminus X$, choose $\phi(v)\in \{1\LL 4p+2\}$, uniformly and independently at random. For $1\le i\le 4p+2$ let $A_i$
be the set of all $v\in A\setminus X$ with $\phi(v) = i$. Thus $X, A_1\LL A_{4p+2}$ are pairwise disjoint sets with union $A$.
We will show that with positive probability, the statement of (2) is satisfied.
For each $I\in \mathcal{I}$ let
$A_I:=\bigcup_{i\in I}A_i$.

There are two cases, depending whether  $|M_1\cupcup M_n|\ge |F|/2$ or not. Suppose first that 
$|M_1\cupcup M_n|\ge |F|/2$. For $I\in \mathcal{I}$ and $1\le j\le n$, we say that $j$ is {\em bad} for $I$
if at most $|M_j|/8$ edges of $M_j$ have both ends in $A_I$. By the first statement of \ref{probability}, since each vertex of $A\setminus X$
belongs to $A_I$ independently with probability $1/2$, and $|M_j|=q$, it follows that the probability that $j$ is bad for $I$ is at most $e^{-q/32}$.
Consequently the expected number of values of $j\in \{1\LL n\}$ such that $j$ is bad for some $I\in \mathcal{I}$ is at most
$$ne^{-q/32}|\mathcal{I}|\le ne^{-q/32}2^{4p+2}\le n/2.$$
%e^{-q/32}2^{2r}\le 1/2
Let $J$ be the set of $j\in \{1\LL n\}$ such that $j$ is not bad for any $I\in \mathcal{I}$. It follows that 
$|J|\ge n/2$ with positive probability. If $|J|\ge n/2$, then 
$$|\bigcup_{j\in J}M_j|\ge |M_1\cupcup M_n|/2\ge |F|/4.$$ 
Moreover, for each $I\in \mathcal{I}$, at least $q/8$ edges of $M_j$ have both ends in $A_I$, for each $j\in J$; and so at least $1/8$
of the edges of $\bigcup_{j\in J}M_j$ have both ends in $A_I$. Consequently, with positive probability at least $|F|/32$ edges of $G[A]$ 
have both ends in $A_I$, and hence in this case the claim is true.

Now we assume that $|M_1\cupcup M_n|\le |F|/2$, and so $|M_0|\ge |F|/2$. 
For each $v\in A\setminus X$, let $d_v$ be the number of neighbours of $v$ in $X$, and let $m=\sum_{v\in A\setminus X}d_v$.
For each $I\in \mathcal{I}$, the probability that $\sum_{v\in A_I}d_v\le m/4$
is at most $e^{-m/(8|X|)}\le e^{-m/(16q)}$, 
by the second statement of \ref{probability}, taking $d=|X|$. But $|F|\ge p|A|\ge 8q^2$ by (1), and so $m\ge |F|/2-2q^2\ge |F|/4\ge 2q^2 $.
Consequently, for each $I\in \mathcal{I}$, the probability that $\sum_{v\in A_I}d_v\le m/4$
is at most $e^{-q/8}$; and hence the probability that $\sum_{v\in A_I}d_v> m/4$ for each $I\in \mathcal{I}$ is at least
$1-2^{4p+2}e^{-q/8}>0$.  We deduce that there is a partition of $A\setminus X$ into $4p+2$ subsets $A_1\LL A_{4p+2}$, such that 
$\sum_{v\in A_I}d_v> m/4$
for each $I\in \mathcal{I}$. But $\sum_{v\in A_I}d_v$ is at most the number of edges that have both ends in $X\cup A_I$. This proves
(2).

\bigskip

Choose the sets $A_1\LL A_{4p+2}$ as in (2), and as before, let
$A_I:=\bigcup_{i\in I}A_i$ for each $I\in \mathcal{I}$. Let $W_0$ be the set of vertices in $V(G)\setminus A$ with
a neighbour in $X$. For each $I\in \mathcal{I}$, let $W_I$ be the set of vertices $v\in V(G)\setminus A$ with no neighbour in
$X\cup A_I$.
From \ref{critical}, every vertex in $V(G)\setminus A$ has at most $2p+1$ neighbours in $A$, and so 
$V(G)\setminus A$ is the union of $W_0$ and the sets $W_I\;(I\in \mathcal{I})$.
Since $G[A]$ has no non-null subgraph with minimum degree at least $2p+1$
%r\ge p
(from the minimality of $A$), it follows that $\chi(A)\le 2p+1$. Also, $\chi(W_0)\le |X|d'\le 2qd'$. Let $I\in \mathcal{I}$.
Thus $G[X\cup A_I]$ has at least $p|A|/32$ edges
(by the choice of $A_1\LL A_{4p+2}$) and at most $|A|$ vertices,
and therefore its denseness is at least $p/32= c$. Since $G[X\cup A_I]$ is anticomplete to $W_I$, we may assume that $\chi(W_I)<c$,
since otherwise the theorem holds. Since $|\mathcal{I}|\le 2^{4p+2}$, it follows that 
$$\chi(G)\le 2p+1+2qd'+2^{4p+2}c\le d,$$
as required. This proves \ref{chimainthm2}.~\bbox

Now we prove \ref{mindegmainthm}, again restated in terms of denseness:

\begin{thm}\label{mindegmainthm2}
For all integers $t, c\ge 1$,
there exists $d\ge 1$, such that if $G$ has denseness at least $d$ and $\tau(G)< t$,
then there are anticomplete subsets $A,B\subseteq V(G)$ where $G[A],G[B]$ both have denseness at least $c$.
\end{thm}
%$p/32\ge c$.
%2t\le p/4 defn of p

%e^{-q/32}2^{3p}\le 1/2

%(2q^2 + 2q^22^{2q}(t-1))\le st defn of s

%$d-4st\ge p$
%$(d/2-(2st + ts^t2^{2t}))/t \ge c$)
%ie d\ge 2tc+ 2(2st + ts^t2^{2t}) defn of d
%d/2-2st\ge 2^{3p}c/2 + (3p+2)/2

%p\ge c
\Proof
This proof shares some ideas with the proof of \ref{chimainthm2}, but has some significant differences. In particular, the 
proof is not by induction on $t$. 
Define $p=\max(32c,4t)$ and let $q$ be an integer with $e^{-q/32}2^{8p+4}\le 1/2$. Choose $s$ with 
$st\ge 2q^2 + 2^{2q+1}q(t-1)$, and choose $d$ with 
$$d>\max\left(p+2st, 2ct+ 2st + ts^t2^{2t},2st+ 2^{8p+4}c + 3p+2\right).$$
We will show that $d$ satisfies the theorem.

Let $G$ be a graph with denseness at least $d$ and $\tau(G)< t$.
Choose vertex-disjoint subsets $R_1\LL R_k$ of $V(G)$ with $k$ maximum, such that for $1\le i\le k$, $R_i$ is a $p$-rock 
of $G\setminus (R_1\cupcup R_{i-1})$ and $|R_i|\le s$. 
\\
\\
(1) {\em $k\le 2t$.}
\\
\\
Suppose that $k\ge 2t$, and let $R_1\cupcup R_{2t} = R$. For $1\le i\le 2t$ let $Z_i$ be the set of all vertices in $V(G)\setminus R_i$ that have no neighbour in
$R_i$. Let $W$ be the set of all $v\in V(G)\setminus R$ that have a neighbour in $R_i$ for at least $t$ values of $i\in \{1\LL 2t\}$.
For each $I\subseteq \{1\LL 2t\}$ with $|I|=t$, and each choice of 
$a_i\in R_i$ for each $i\in I$, there are fewer than $t$ vertices in $V(G)\setminus R$ adjacent to $a_i$ for each $i\in I$,
since $\tau(G)<t$.
For each $I$ there are at most $s^t$ choices of the vertices $a_i\;(i\in I)$, and so there are at most $ts^t$ vertices in $V(G)\setminus R$ 
with a neighbour in $R_i$ for each $i\in I$. Since there are at most $2^{2t}$ choices of $I$, it follows that $|W|\le ts^t2^{2t}$.
Thus $|R\cup W|\le 2st + ts^t2^{2t}$, and so at most $(2st + ts^t2^{2t})|G|$ edges have an end in $R\cup W$. Since $G$ has at least 
$d|G|$ edges, there are at least $(d-(2st + ts^t2^{2t}))|G|$ edges with neither end in $R\cup W$. For every such edge, say $uv$, since $u$ has a neighbour
in at most $t-1$ of $R_1\LL R_{2t}$, and the same for $v$, there exists $i\in \{1\LL 2t\}$ such that neither of $u,v$ has a neighbour in $R_i$, that is, $u,v\in
Z_i$.
Consequently there exists $i\in \{1\LL 2t\}$ such that at least $(d-(2st + ts^t2^{2t}))|G|/(2t)$ edges $uv$ of $G$ have both ends in $Z_i$. It follows that
$G[Z_i]$ has denseness at least $(d-(2st + ts^t2^{2t}))/(2t)\ge c$, and it is anticomplete to the rock $R_i$, and so the theorem holds.
%$(d/2-(2st + ts^t2^{2t}))/t \ge c$
%ie d\ge 2tc+ 2(2st + ts^t2^{2t})
This proves (1).

\bigskip
Let $R=R_1\cupcup R_k$. Thus $|R|\le 2st$ by (1). Consequently at most $2st|G|$ edges of $G$ have an end in $R$, and so the graph $G\setminus R$ has at least
$(d-2st)|G|$ edges.
%$d-4st\ge p$
Since $d-2st\ge p$, there is a rock $A$ of $G\setminus R$. From the maximality of $k$, $|A|> s$.

From \ref{getmatching}, there is a partition of $E(G[A])$ into sets $M_0\LL M_n$ for some $n\ge 0$,
such that
\begin{itemize}
\item there is a subset $X\subseteq V(A)$ with $|X|\le 2q-2$
such that every edge in $M_0$ is incident with a vertex in $X$; and
\item $M_1\LL M_n$ are all matchings, each with cardinality $q$.
\end{itemize}

\bigskip
\noindent
(2) {\em $|M_0|\le 2t|A|\le p|A|/2$, and hence $M_1\cupcup M_n$ has cardinality at least $p|A|/2$.}
%2t\le p/4
\\
\\
There are at most $2q^2$ edges in $E(G[A])$ with both ends in $X$, since $|X|\le 2q$. We need to count the number with exactly one end in $X$.
For each subset $Y$ of $X$ with $|Y|=t$, there are at most $t-1$ vertices adjacent to each vertex in $Y$, and so there are at most $2^{2q}(t-1)$ vertices in
$A\setminus X$ with at least $t$ neighbours in $X$. Hence there are at most $2^{2q}(t-1)|X|\le 2^{2q+1}q(t-1)$ edges $uv$ of $G[A]$ with $u\in X$ and $v\in A\setminus X$
such that $v$ has at least $t$ neighbours in $X$. But there are at most $(t-1)|A|$ edges $uv$ of $G[A]$ with $u\in X$ and $v\in A\setminus X$
such that $v$ has fewer than $t$ neighbours in $X$; so altogether there are at most
$$2q^2 + 2^{2q+1}q(t-1) + (t-1)|A|\le \left((2q^2 + 2^{2q+1}q(t-1))/s + (t-1)\right)|A|\le 2t|A|\le p|A|/2$$
edges of $G[A]$ with an end in $X$, since $|A|\ge s$. This proves the first statement of (2). The second follows since $|E(G[A])|\ge p|A|$. This proves (2).
%(2q^2 + 2q^22^{2q}(t-1))\le st

\bigskip
Let $\mathcal{I}$ be the set of all subsets of $\{1\LL 8p+4\}$ with cardinality $4p+2$.
\\
\\
(3) {\em There is a partition of $A\setminus X$ into $8p+4$ subsets $A_1\LL A_{8p+4}$, such that for each $I\in \mathcal{I}$ there are at
least $p|A|/32$ edges of $G[A]$ that have both ends in $\bigcup_{i\in I}A_i$.}
\\
\\
For each $v\in A\setminus X$, choose $\phi(v)\in \{1\LL 8p+4\}$, uniformly and independently at random. For $1\le i\le 8p+4$ let $A_i$
be the set of all $v\in A\setminus X$ with $\phi(v) = i$. Thus $X, A_1\LL A_{8p+4}$ are pairwise disjoint sets with union $A$.
We will show that with positive probability, the statement of (3) is satisfied.
For each $I\in \mathcal{I}$ let
$A_I:=\bigcup_{i\in I}A_i$.

For $I\in \mathcal{I}$ and $1\le j\le n$, we say that $j$ is {\em bad} for $I$
if at most $q/8$ edges of $M_j$ have both ends in $A_I$. By the first statement of \ref{probability}, since each vertex of $A\setminus X$
belongs to $A_I$ independently with probability $1/2$, it follows that the probability that $j$ is bad for $I$ is at most $e^{-q/32}$.
Consequently the expected number of values of $j\in \{1\LL n\}$ such that $j$ is bad for some $I\in \mathcal{I}$ is at most
$$ne^{-q/32}|\mathcal{I}|\le ne^{-q/32}2^{8p+4}\le n/2.$$
%e^{-q/32}2^{3p}\le 1/2
Let $J$ be the set of $j\in \{1\LL n\}$ such that $j$ is not bad for any $I\in \mathcal{I}$. It follows that
$|J|\ge n/2$ with positive probability. Moreover, if $|J|\ge n/2$, then 
$$|\bigcup_{j\in J}M_j|\ge |M_1\cupcup M_n|/2\ge p|A|/4$$ 
by (2).
But for each $I\in \mathcal{I}$, at least $q/8$ edges of $M_j$ have both ends in $A_I$, for each $j\in J$; and so at least $1/8$
of the edges of $\bigcup_{j\in J}M_j$ have both ends in $A_I$. Consequently, with positive probability at least $p|A|/32$ edges of $G[A]$ 
have both ends in $A_I$. This proves (3).

\bigskip

Choose $A_1\LL A_{8p+4}$ as in (3), and as before, let
$A_I:=\bigcup_{i\in I}A_i$ for each $I\in \mathcal{I}$.
For each $I\in \mathcal{I}$, let $W_i$ be the set of vertices in $V(G\setminus (A\cup R))$ with no neighbour in $A_I$. Since for every edge 
$uv$ of $G\setminus R$ with $u,v\notin A$, $u$ has a neighbour in $A_i$ for at most $2p+1$ values of $i\in \{1\LL 8p+6\}$ by \ref{critical}, 
and the same for $v$, it follows
that there exists $I\in \mathcal{I}$ with $u,v\in W_i$. But, since $G[A_I]$ has denseness at least $p/32\ge c$ by (3), 
%$p/32\ge c$
and is anticomplete to $W_I$, we may assume that $G[W_I]$ has denseness less than $c$, and so there are at most $c|G|$
edges of $G\setminus R$ with both ends in $A_I$. We will show that this leads to a contradiction.
Since there are only at most $2^{8p+4}$ choices of $I$, there are at most 
$2^{8p+4}c|G|$ edges of $G\setminus R$ with neither end in $A$. But there are at most $(2p+1)|G|$ edges with one end in $A$ and the other in 
$V(G)\setminus (A\cup R)$, since every vertex in $V(G)\setminus (A\cup R)$ has at most $2p+1$ neighbours in $A$ by \ref{critical}. Also, from the 
minimality of $A$ (in the definition of a rock), if we delete a vertex of $A$, the remainder induces a graph with fewer than
$p(|A|-1)$ edges, and so $G[A]$ has fewer than
$$p(|A|-1) + |A|\le (p+1)|A| \le (p+1)|G|$$
edges. Altogether, then, $G\setminus R$ has fewer than
$$2^{8p+4}c|G| + (2p+1)|G|+ (p+1)|G|<(d-2st)|G|$$
%d/2-2st\ge 2^{3p}c/2 + (3p+2)/2
edges. But we already saw that $G\setminus R$ has at least $(d-2st)|G|$ edges, a contradiction. This proves \ref{mindegmainthm2}.~\bbox

\section{Tournaments}
There is an interesting extension of \ref{erdosconj} to tournaments. 
If $G$ is a tournament, a subset $X\subseteq V(G)$ is {\em acyclic} if it has no directed cycle; and $\chi(G)$
is the minimum $k$ such that $V(G)$ is the union of $k$ acyclic subsets. 
Again, we write $\chi(A)$ for $\chi(G[A])$
when $A\subseteq V(G)$. 
If $A,B\subseteq V(G)$ are disjoint, we say $A$ 
is {\em complete} to $B$ if every vertex in $B$ is adjacent from every vertex in $A$. 
We have not been able to find a counterexample 
to the 
following strengthening of \ref{erdosconj}.
\begin{thm}\label{tourconj2}
{\bf Conjecture: }For all integers $c\ge 1$ there exists $d\ge 1$ such that if $G$ is a tournament and $\chi(G)\ge d$, there are
disjoint $A,B\subseteq V(G)$, with $A$ complete to $B$, and both inducing tournaments with chromatic number at least $c$.
\end{thm}
We will discuss this further in another paper~\cite{tourprobs}, where we will prove that it implies \ref{erdosconj}, and prove the following two results:
\begin{thm}\label{triangleform}
For all $c\ge 1$ there exists $d\ge 1$ such that if $G$ is a tournament with $\chi(G)\ge d$, then there exist disjoint
$A,B\subseteq V(G)$ with $A$ complete to $B$, where $A$ is a cyclic triangle and $\chi(B)\ge c$.
\end{thm}
(A {\em cyclic triangle} is a three-vertex set inducing a directed cycle.)
The second result concerns domination number. A tournament $G$ has {\em domination number} $k$ if $k$ is minimum such that for some set
$X\subseteq V(G)$ with $|X|=k$, every vertex in $V(G)\setminus X$ is adjacent from some vertex in $X$.
\begin{thm}\label{bigdom}
For every integer $c\ge 1$, there exists $d\ge 1$ such that if $G$ is a tournament with domination number at least $d$,
then there are
disjoint $A,B\subseteq V(G)$, such that $A$ is complete to $B$
and  $\chi(A), \chi(B)\ge c$.
\end{thm}

\end{document}